\documentclass[12pt,a4paper]{amsart}
\usepackage{amsfonts}

\usepackage{amscd,amssymb}
\usepackage[pdftex]{graphicx}
\usepackage{times}
\usepackage[T1]{fontenc}

\theoremstyle{plain}
\newtheorem{theorem}{Theorem}[section]
\newtheorem{corollary}[theorem]{Corollary}
\newtheorem{lemma}[theorem]{Lemma}

\theoremstyle{definition}
\newtheorem{definition}{Definition}[section]

\theoremstyle{example}
\newtheorem{example}{Example}[section]

\numberwithin{equation}{section}

\numberwithin{table}{section}

\usepackage[
    pdftex,
    pdfpagemode=None,
    bookmarksopen=true,
    baseurl={http://ajmaa.org/}
    pdfstartview={FitH}
]{hyperref}

\theoremstyle{plain}
\newtheorem*{acknowledgement}{Acknowledgements}

\numberwithin{equation}{section}

\newcommand{\eps}{\mathbb{\varepsilon}}

\newcommand{\R}{\mathbb{R}}
\newcommand{\rar}{\mbox{$\rightarrow$}}

\newcommand{\ti}{\tilde}
\newcommand{\T}{\mathbb{T}}

\newcommand{\Z}{\mathbb{Z}}

\begin{document}

\title[The Riemann-Stieltjes integral on time scales]{The Riemann-Stieltjes integral on time scales\footnote{Accepted
(March 6, 2009) for publication at
\emph{The Australian Journal of Mathematical Analysis and
Applications} (AJMAA).}}


\author[D. Mozyrska]{Dorota Mozyrska}
\address{Faculty of Computer Science, Bia{\l}ystok Technical University, 15-351 Bia\l ystok, Poland}
\email{\href{mailto: Dorota Mozyrska <admoz@w.tkb.pl>}{admoz@w.tkb.pl}}


\author[E. Paw\l uszewicz]{Ewa Paw\l uszewicz}
\address{Faculty of Computer Science, Bia{\l}ystok Technical University, 15-351 Bia\l ystok, Poland}
\curraddr{Department of Mathematics, University of Aveiro,
		 3810-193 Aveiro, Portugal}
\email{\href{mailto: Ewa Pawluszewicz <epaw@pb.edu.pl>}{epaw@pb.edu.pl},
\href{mailto: Ewa Pawluszewicz <ewa@ua.pt>}{ewa@ua.pt}}


\author[D. F. M. Torres]{Delfim F. M. Torres}
\address{Department of Mathematics, University of Aveiro,
		 3810-193 Aveiro, Portugal}
\email{\href{mailto: Delfim F. M. Torres <delfim@ua.pt>}{delfim@ua.pt}}


\date{}




\begin{abstract}
We study the process of integration on time scales
in the sense of Riemann-Stieltjes. Analogues of the
classical properties are proved for a generic time scale,
and examples are given.
\end{abstract}

\maketitle


\section{Introduction}

The development of integration on time scales is a recent
but already well studied subject. Available integration notions
on time scales include the Riemann delta integral \cite{Bh1,G:K:2002,G:K:2004},
the Riemann nabla integral \cite{G:K:2002}, the Riemann diamond-alpha integral
\cite{M:T:DSA}, the Lebesgue delta and nabla integrals \cite{Bh1,G:2003},
and the more general Henstock-Kurzweil delta and nabla integrals \cite{P:2006,Thomson}.
Other studies on time scales are dedicated to improper integrals \cite{B:G:2003}
and to multiple integration \cite{B:G:2005,B:G:2006}. Surprisingly enough, the
Stieltjes integral has not received attention in the literature of time scales.

In this paper we study the process of Stieltjes integration
on time scales, both in nabla and delta sense.
We trust that such integrals will find interesting applications in the study of dynamic equations on time scales, enabling to study more general situations than those treated before (\textrm{cf.} Section~\ref{sec:ConcFuture}).
Because delta and nabla theories are similar,
we avoid repetition by following the approach promoted by Bartosiewicz and Piotrowska \cite{BPi}: the box symbol $\Box$ is here used to represent the delta operator $\Delta$ as well as the nabla operator $\nabla$.

It is assumed that the reader is familiar with the time
scale calculus and the notations for delta and nabla differentiation.
For an introduction to time scales the reader is referred
to the book by Bohner and Peterson \cite{Bh}.
Throughout the paper $\T$ denotes a time scale.
Let $a,b\in\T$ and $a<b$. We distinguish $[a,b]$ as a real interval and we define $I=[a,b]_{\T}:=[a,b]\cap\T$.
In this sense, $[a,b]=[a,b]_{\R}$. Along the text
$I$ is a nonempty, closed, and bounded interval consisting
of points from a time scale $\T$.
Moreover, if $I=[a,b]_{\T}$, then we define
$I_{\Delta}:=[a,\rho(b)]_{\T}$ and $I_{\nabla}:=[\sigma(a), b]_{\T}$. By $I_{\Box}$ we denote one of them, where $\Box$
means either $\nabla$ or $\Delta$. Similarly, we use "$\Box$"
as  a common notation for the two kinds of derivatives on time scales: one can read $f^{\Box}$ either as $f^{\Delta}$
or as $f^{\nabla}$.


\section{The Riemann-Stieltjes $\Box$--integral}
\label{sec:RS}

Let $\T$ be a time scale, $a,b\in\T$, $a<b$,
and $I=[a,b]_{\T}$. A partition of $I$ is any finite ordered subset
\[
P=\{t_0, t_1,\ldots, t_n\}\subset [a,b]_{\T}, \quad
\text{where} \  a=t_0<t_1<\ldots<t_n=b\, .
\]
Let  $g$ be a real-valued increasing function on $I$. Each partition
$P=\{t_0, t_1,\ldots, t_n\}$ of $I$ decomposes $I$ into
subintervals $I_{\Box j}=[t_{j-1},t_j]_{\Box}$, $j=1,2,\ldots, n$, such
that $I_{\Box j}\cap I_{\Box k}=\emptyset$ for any $k\neq j$.
By $\Delta t_j=t_j-t_{j-1}$ we denote the length of the $j$th
subinterval in the partition $P$; by $\mathcal{P}(I)$ the
set of all partitions of $I$.

Let $P_m$, $P_n \in\mathcal{ P}(I)$. If $P_m\subset P_n$ we call $P_n$
a refinement of $P_m$. If $P_m$, $P_n$ are independently chosen, then
the partition $P_m\cup P_n$ is a common refinement of $P_m$ and
$P_n$.

Let us now consider a strictly increasing real-valued function $g$
on the interval $I$. Then, for the partition $P$ of $I$ we define
\[g(P)=\{g(a)=g(t_0), g(t_1), \ldots, g(t_{n-1}),g(t_n)=g(b)\}\subset g(I)\]
and $\Delta g_j=g(t_j)-g(t_{j-1})$. We note that $\Delta g_j$ is
positive and $\sum\limits_{j=1}^n\Delta g_j=g(b)-g(a)$. Moreover,
$g(P)$ is a partition of $[g(a),g(b)]_{\R}$. In what follows, for
the particular case $g(t)=t$ we obtain the Riemann sums for delta
integrals studied in \cite{Bh1}. We note that for a general $g$ the image $g(I)$ is not necessarily an interval in the classical sense, even for rd-continuous functions $g$, because our interval $I$ may contain scattered points. From now on let $g$ be always a strictly increasing real function on the considered interval $I=[a,b]_{\T}$.

\begin{lemma}
\label{lem:2}
Let $I=[a,b]_{\T}$ be a closed (bounded) interval in $\T$ and let $g$ be continuous on $I$. For every $\delta>0$
there is a partition $P_{\delta}=\{t_0,t_1,\ldots, t_n\}\in
\mathcal{P}(I)$ such that for each $j\in\{1,2,\ldots,n\}$ one has:
\[\Delta g_j=g(t_j)-g(t_{j-1})\leq \delta \quad \mbox{or} \quad
\Delta g_j>\delta \ \wedge \ \rho(t_j)=t_{j-1}.\]
\end{lemma}
\begin{proof}
Let $I=[a,b]_{\T}$ be closed (bounded) interval in $\T$. Firstly, let
us observe that if for $t\in I$ $g(t)\in [g(a),g(b)]_{\R}$, then
$a\leq t\leq b$ because $g$ is an increasing function on $I$. We
define two families of sets,
$B_j=(g(t_{j-1}),g(t_{j-1})+\delta]_{\R}\cap [g(a),g(b)]_{\R}$ and
$A_j=\{t\in I: g(t)\in B_j\}$. Inductively, we
construct the partition taking $t_0=a$ and
\[t_j=\left\{\begin{array}{lr} \sup A_j,  & \mbox{if} \ A_j\neq \emptyset \, ,\\
		\sigma(t_{j-1}), & \mbox{if} \ A_j= \emptyset \, .
\end{array}\right. \]
Then, we get $t_0<t_1<\ldots<t_n$, where $t_n=b$. It follows that
if $t_{j-1}<t_j$ and $t_j=\sigma(t_{j-1})$,
then $\rho(t_j)=t_{j-1}$ for any time scale $\T$.
\end{proof}

Let $f$ be a real-valued and bounded  function on the interval
$I$. Let us take a partition $P=\{t_0, t_1,\ldots, t_n\}$ of $I$.
Denote $I_{\Box j}=[t_{j-1},t_j]_{\Box}$, $j=1,2,\ldots,n$, and
\[m_{\Box j}=\inf\limits_{t\in I_{\Box j}} f(t), \quad  M_{\Box j}=\sup\limits_{t\in I_{\Box j}} f(t)
\, . \]
The \emph{upper Darboux--Stieltjes $\Box$--sum of $f$ with respect to the
partition $P$}, denoted by $U_{\Box}(P,f,g)$, is defined by
\[U_{\Box}(P, f, g)=\sum_{j=1}^nM_{\Box j}\Delta g_j \, ,\]
while the \emph{lower Darboux--Stieltjes $\Box$--sum of $f$ with respect to the
partition $P$}, denoted by $L_{\Box}(P,f,g)$, is defined by
\[L_{\Box}(P, f, g)=\sum_{j=1}^n m_{\Box j}\Delta g_j\,.\]

\begin{example}
\label{ex1}
Let $\T=\overline{q^{\Z}}$, $q>1$, $f(t)=t$,
$g(t)=t^2$, and $I=[0,1]_{\T}$. For the partition
$P=\{0,q^{-n+1}, \ldots, q^{-1},1\}$, where
$t_0=0<q^{-n+1}<\ldots<q^{j-n}<\ldots<q^{-1}<1=t_n$ with
$t_j=q^{j-n}$ for $j=1,\ldots, n$, we have $\Delta
g_j=t_j^2-t_{j-1}^2=q^{2(j-n-1)}(q^2-1)$  for $j= 2,\ldots,n$,
while $\Delta g_1=t_1^2-0=q^{2-2n}$.
Let us read $\Box=\Delta$. In this case
$M_{\Delta j}=\rho(t_j)$ and $m_{\Delta j}=t_{j-1}$. For our partition we have that
$I_{\Delta j}=[t_{j-1},\rho(t_j)]_{\T}=t_{j-1}$, $j= 2,\ldots,n$. Hence,
$m_{\Delta j}=M_{\Delta j}=t_{j-1}$ for $f(t)=t$, $j=2,\ldots,n$, and
$m_{\Delta 1}=0$, $M_{\Delta 1}=\rho(t_1)=q^{-n}$. The lower and upper
Darboux--Stieltjes $\Delta$--sums are, respectively,
\[L_{\Delta}(P,f,g)=\sum\limits_{j=2}^n t_{j-1}\Delta g_j
=\frac{q+1}{q^2+q+1}(1-q^{3(1-n)})\,  ,\]
\[U_{\Delta}(P,f,g)=q^{2-3n}+\sum\limits_{j=2}^n t_{j-1}\Delta g_j
=q^{2-3n}+L(P,f,g)=\frac{q+1+q^{2-3n}}{q^2+q+1}\, .\]
Consider now $\Box=\nabla$. In this case
$M_{\nabla j}=t_j$ and $m_{\nabla j}=\sigma(t_{j-1})$. For our partition we have that
$I_{\nabla j}=[\sigma(t_{j-1}),t_j]_{\T}=t_{j}$, $j=2,\ldots,n$. Hence,
$m_{\nabla j}=M_{\nabla j}=t_{j}$ for $f(t)=t$, $j=2,\ldots,n$, and
$m_{\nabla 1}=\sigma(0)=0$, $M_{\nabla 1}=t_1=q^{1-n}$. The lower and upper Darboux--Stieltjes
$\nabla$--sums are, respectively,
\[L_{\nabla}(P,f,g)=\sum\limits_{j=2}^n t_j \Delta g_j
=q\frac{q+1}{q^2+q+1}(1-q^{3(1-n)})\,  ,\]
\[U_{\nabla}(P,f,g)=q^{3-3n}+\sum\limits_{j=2}^n t_{j}\Delta g_j
=q^{3-3n}+L_{\nabla}(P,f,g)=q\frac{q+1+q^{2-3n}}{q^2+q+1}\, .\]
\end{example}

For computing the value of a Riemann integral on time scales one
uses the fact that $\Delta t_j\geq 0$. Since we assume that
$g$ is an increasing function, then $\Delta g_j\geq 0$. As we
shall see below, also some other properties of the Riemann-integral
(for function $g(t)=t$) are preserved for an arbitrary increasing
function $g$.

\begin{theorem}
\label{lem:1}
Suppose that $f$ is a bounded function on
$I=[a,b]_{\T}$, $a$, $b\in\T$.
Let $P\in \mathcal{P}(I)$ and $m\leq f(t)\leq M$ for all $t\in I$. Then,
\[m(g(b)-g(a))\leq L_{\Box}(P,f,g)\leq U_{\Box}(P,f,g)\leq M(g(b)-g(a))\]
and
\[L_{\Box}(P,f,g)\leq L_{\Box}(Q,f,g)\leq U_{\Box}(Q,f,g)\leq U_{\Box}(P,f,g)\]
for any refinement $Q$ of $P$.
\end{theorem}
\begin{proof}
For $\Box=\Delta$ the proof is an immediate consequence of \cite[Lemma~5.2]{Bh1}. For $\Box=\nabla$ the proof is done
in similar steps.
\end{proof}

\begin{definition}
\label{RSdef}
Let $I=[a,b]_{\T}$, where $a,b\in\T$.
The upper Darboux-Stieltjes $\Box$--integral from $a$ to $b$
with respect to function $g$ is defined by
\[\overline{\int_a^b}f(t)\Box g(t)=\inf_{P\in \mathcal{P}(I)} U_{\Box}(P,f,g)\, ;\]
the lower Darboux-Stieltjes $\Box$--integral from $a$ to
$b$ with respect to function $g$ is defined by
\[\underline{\int_a^b}f(t)\Box g(t)=\sup_{P\in \mathcal{P}(I)} L_{\Box}(P,f,g)\, .\]
If  $\overline{\int_a^b}f(t)\Box
g(t)=\underline{\int_a^b}f(t)\Box g(t)$, then we say that $f$ is
$\Box$--integrable with respect to $g$ on $I$, and the common
value of the integrals, denoted by $\int_a^bf(t)\Box
g(t)=\int_a^bf\Box g$, is called the Riemann-Stieltjes (or
simply Stieltjes) $\Box$--integral of $f$ with respect to $g$ on
$I$.

The set of all functions that
are $\Box$--integrable with respect to $g$ in the Riemann-Stieltjes
(also Darboux-Stieltjes) sense will be denoted by
$\mathcal{R}_{\Box}(g,I)$.
\end{definition}

\begin{theorem}
\label{th.7}
If  $L_{\Box}(P,f,g)=U_{\Box}(P,f,g)$ for some $P \in \mathcal{P}(I)$,
then function $f$ is Riemann-Stieltjes $\Box$--integrable
and \[\int_a^bf\Box g=L_{\Box}(P,f,g)=U_{\Box}(P,f,g) \, .\]
\end{theorem}
\begin{proof}
Follows immediately from the definition of Riemann-Stieltjes
$\Box$--integral and Corollary~\ref{lem:1}.
\end{proof}

\begin{example}
Let $\T=\overline{q^{\Z}}$, $q>1$,
and let us continue Example~\ref{ex1}.
Consider functions $f(t)=t$ and $g(t)=t^2$ on the interval $I=[0,1]_{\T}$.
For the partition $P=\{0,q^{-n+1}, \ldots, q^{-1},1\}$, where
$t_0=0<q^{-n+1}<\ldots<q^{j-n}<\ldots<q^{-1}<1=t_n$ with $t_j=q^{j-n}$ for $j=1,\ldots, n$,
we have
\[\overline{\int_a^b}f(t)\Delta g(t)=\inf\limits_{P\in \mathcal{P}(I)} U_{\Delta}(P,f,g)
=\lim\limits_{n\rightarrow\infty}\frac{q+1+q^{2-3n}}{q^2+q+1}=\frac{q+1}{q^2+q+1}\]
and
\[\underline{\int_a^b}f(t)\Delta g(t)=\sup\limits_{P\in \mathcal{P}(I)} L_{\Delta}(P,f,g)
=\lim\limits_{n\rightarrow\infty}\frac{q+1}{q^2+q+1}(1-q^{3(1-n)})
=\overline{\int_a^b}f(t)\Delta g(t)\, .\]
Consequently, $\int_a^b f(t)\Delta g(t)=\frac{q+1}{q^2+q+1}$.
We also have that
\[\overline{\int_a^b}f(t)\nabla g(t)=\inf\limits_{P\in \mathcal{P}(I)} U_{\nabla}(P,f,g)=\frac{q^2+q}{q^2+q+1}\]
and
\[\underline{\int_a^b}f(t)\nabla g(t)=\sup\limits_{P\in \mathcal{P}(I)} L_{\nabla}(P,f,g)
=\lim\limits_{n\rightarrow\infty}\frac{q^2+q}{q^2+q+1}(1-q^{3(1-n)})
=\overline{\int_a^b}f(t)\nabla g(t)\, .\]
Hence, $\int_a^b f(t)\nabla g(t)=\frac{q^2+q}{q^2+q+1}=q\int_a^b f(t)\Delta g(t)$.
\end{example}

\begin{theorem}
Let $f$ be a bounded function on $I=[a,b]_{\T}$, $a$, $b\in \T$,
$m\leq f(t)\leq M$ for all $t\in I$,
and $g$ be a function defined and monotonically increasing on $I$. Then,
\[m(g(b)-g(a))\leq \underline{\int_a^b}f(t)\Box g(t)\leq
\overline{\int_a^b}f(t)\Box g(t)\leq M(g(b)-g(a))\, .\]
If $f\in \mathcal{R}_{\Box}(g,I)$, then
\[m(g(b)-g(a))\leq \int_a^bf(t)\Box g(t)\leq  M(g(b)-g(a))\, .\]
\end{theorem}
\begin{proof}
The definition of upper and lower Darboux-Stieltjes
$\Box$--integral implies that
$\underline{\int_a^b}f(t)\Box g(t) \leq
\overline{\int_a^b}f(t)\Box g(t)$. Thus,
\[\overline{\int_a^b}f(t)\Box g(t)=\inf_{P\in \mathcal{P}(I)}
U_{\Box}(P,f,g)\leq U_{\Box}(P,f,g) \leq \sum\limits_{j=1}^nM \Delta
 g_j=M(g(b)-g(a)) \, .\]
Similarly, $m(g(b)-g(a)) \leq \underline{\int_a^b}f(t)\Box g(t)$,
and the proof is done.
\end{proof}

\begin{theorem}(Integrability criterion)
\label{dor:criterion}
Let $f$ be a bounded function on $I=[a,b]_{\T},a,b\in \T$. Then,
$f\in\mathcal{R}_{\Box}(g,I)$ if and only if for every $\eps>0$
there exists a partition
$P\in \mathcal{P}(I)$ such that
\begin{equation}
\label{dor:criterion1}
U_{\Box}(P,f,g)-L_{\Box}(P,f,g)<\eps\, .
\end{equation}
\end{theorem}
\begin{proof}
Suppose that $f\in \mathcal{R}_{\Box}(g,I)$ and let $\eps>0$. Because
$$\int_a^b f\Box g=\inf\limits_{P\in
\mathcal{P}(I)}U_{\Box}(P,f,g)=\sup\limits_{P\in \mathcal{P}(I)} L_{\Box}(P,f,g)$$
there exist $P_1$ and $P_2\in \mathcal{P}(I)$ such that
\[\int_a^b f\Box g<U_{\Box}(P_1,f,g)< \int_a^b f\Box g+\frac{\eps}{2}\]
and
\[\int_a^b f\Box g-\frac{\eps}{2}<L_{\Box}(P_2,f,g)< \int_a^b f\Box g\, .\]
Hence, $U_{\Box}(P_1,f,g)-\int_a^b f(t)\Box g(t)<\frac{\eps}{2}$ and
$\int_a^b f(t)\Box g(t)-L_{\Box}(P_2,f,g)<\frac{\eps}{2}$. Let $P$ be
a common refinement of $P_1$ and $P_2$. Then, Theorem~\ref{lem:1}
implies that $U_{\Box}(P,f,g)-\int_a^b f(t)\Box g(t)<\frac{\eps}{2}$ and
$\int_a^b f(t)\Box g(t)-L_{\Box}(P,f,g)<\frac{\eps}{2}$. Thus,
$U_{\Box}(P,f,g)-L_{\Box}(P,f,g)<\eps$.

Conversely, let $f$ be a bounded function on $I$ and $g$ be an
increasing function on $I$. Suppose that for $\eps>0$ there exists
$P\in \mathcal{P}(I)$ such that $U_{\Box}(P,f,g)-L_{\Box}(P,f,g)<\eps$. The definition
of Riemann-Stieltjes $\Box$--integral and Theorem~\ref{lem:1}
imply that
\[L_{\Box}(P,f,g)\leq \underline{\int_a^b} f\Box g
\leq \overline{\int_a^b} f\Box g\leq U_{\Box}(P,f,g)\, .\]
It is obvious that $0\leq \overline{\int_a^b} f\Box
g-\underline{\int_a^b} f\Box g<\eps$. Since $\eps$ is
arbitrary, then $\overline{\int_a^b} f(t)\Box
g(t)$ $=\underline{\int_a^b} f\Box g$ and
$f\in\mathcal{R}_{\Box}(g,I)$.
\end{proof}

\begin{theorem}
Suppose that $P_{\delta}$ is a partition as stated in
Lemma~\ref{lem:2}. A bounded function $f$ on $I=[a,b]_{\T}$ is
integrable if and only if for each $\eps>0$ there exists  $\delta>0$
such that
\begin{equation*}
P_{\delta}\in \mathcal{P}(I)
\Rightarrow
U_{\Box}(P_{\delta},f,g)- L_{\Box}(P_{\delta},f,g)<\eps \, .
\end{equation*}
\end{theorem}
\begin{proof}
Part "$\Leftarrow$" follows immediately from the integrability
criterion (Theorem \ref{dor:criterion}).

For the proof of the part "$\Rightarrow$" suppose that $f \in
\mathcal{R}_{\Box}(g,I)$. Let $\varepsilon>0$. Then, there exist
partitions $\mathcal{P}_{\delta}^1$ and $\mathcal{P}_{\delta}^2$
such that
\begin{equation}
\label{1}
 U(\mathcal{P}_{\delta}^1,f,g)<\overline{\int_a^b} f\Box g+\frac{\varepsilon}{2}
\end{equation}
and
\begin{equation}
\label{2}
 \underline{\int_a^b}f\Box g-\frac{\varepsilon}{2}<L(\mathcal{P}_{\delta}^2,f,g).
\end{equation}
Let
$\mathcal{P}_{\delta}=\mathcal{P}_{\delta}^1\cup\mathcal{P}_{\delta}^2$.
Theorem~\ref{lem:1} and inequalities (\ref{1}) and (\ref{2}) imply that
\[\underline{\int_a^b}f\Box g-\frac{\varepsilon}{2}<L(\mathcal{P}_{\delta},f,g) \leq
  U(\mathcal{P}_{\delta},f,g)<\underline{\int_a^b}f\Box g+\frac{\varepsilon}{2}\]
Because $f \in \mathcal{R}_{\Box}(g,I)$, then
$\underline{\int_a^b}f\Box g=\underline{\int_a^b}f\Box g$. Hence,
$U(\mathcal{P}_{\delta},f,g)-L(\mathcal{P}_{\delta},f,g)<\varepsilon$.
\end{proof}

Lemma~\ref{lem:1} and the properties of the Riemann delta (nabla) integral imply the following:
\begin{theorem}
Let $I=[a,b]_{\T}$, $a$, $b\in\T$. Then, the condition $f\in
\mathcal{R}_{\Box}(g,I)$ is equivalent to each one of the following items:
\begin{enumerate}
  \item [i)]  $f$ is a monotonic function on $I$;

  \item [ii)] $f$ is a continuous function on $I$;

  \item [iii)] $f$ is regulated on $I$;

  \item [iv)] $f$ is bounded and has a finite number of discontinuity points on $I$.
\end{enumerate}
\end{theorem}

\begin{theorem}
\label{dor:3}
Let $f$ be bounded on $I=[a,b]_{\T}$, $a<b$. Then,
\begin{enumerate}
 \item[i)] If there exists an $\eps>0$ and a partition
 $P^*\in \mathcal{P}(I)$ such that the inequality (\ref{dor:criterion1})
 in Theorem~\ref{dor:criterion} is satisfied, then (\ref{dor:criterion1})
 is satisfied for every refinement $P$ of $P^*$.

 \item[ii)] If inequality (\ref{dor:criterion1})
 is satisfied for the partition $P$ given by
 $t_0=a<t_1<\ldots <t_n=b$, and $\xi_j$,
 $\tau_j\in I_{\Box j}=[t_{j-1},t_j]_{\Box}$ for $j=1, 2,\ldots, n$, then
  \begin{equation*}
   \sum_{j=1}^n\left|f(\xi_j)-f(\tau_j)\right|\Delta g_j<\eps\, .
  \end{equation*}

 \item[iii)] If $f\in \mathcal{R}_{\Box}(g,I)$ and inequality
 (\ref{dor:criterion1}) is satisfied for the partition $P$
 given by $t_0=a<t_1<\ldots <t_n=b$ and $\xi_j\in I_{\Box j}=[t_{j-1},t_j]_{\Box}$ for $j=1, 2, \ldots, n$, then
  \begin{equation*}
   \left|\sum_{j=1}^nf(\xi_j)\Delta g_j-\int_a^b f(t)\Box g(t)\right|<\eps\, .
  \end{equation*}
\end{enumerate}
\end{theorem}
\begin{proof}
Part i) follows from Theorem~\ref{lem:1}. Let us now prove part
ii). Let $f$ be a real-valued and bounded function on the interval
$I$. Let us take a partition $P=\{t_0, t_1,\ldots, t_n\}$ of $I$.
 Let $\xi_j$, $\tau_j\in
I_{\Box j}=[t_{j-1},t_j]_{\Box}$ for $j=1, 2,\ldots, n$. Then,
$M_{\Box j}-f(\xi_j)\geq 0$, $f(\tau_j)-m_{\Box j}\geq 0$, and
\[\sum_{j=1}^n\left|f(\xi_j)-f(\tau_j)\right|\Delta g_j
\leq \sum_{j=1}^n\left|\left(f(\xi_j)-f(\tau_j)\right)
+ \left(M_{\Box j}-f(\xi_j)\right)+\left(f(\tau_j)-m_{\Box j}\right)\right|\Delta g_j \, .\]
It follows that
\[\sum_{j=1}^n\left|f(\xi_j)-f(\tau_j)\right|\Delta g_j
\leq \sum_{j=1}^n\left|M_{\Box j}-m_{\Box j}\right|\Delta g_j=U_{\Box}(P,f,g)-L_{\Box}(P,f,g)< \eps \, .\]

Finally, let us prove part iii). Let the inequality
(\ref{dor:criterion1}) be satisfied for the partition $P$ given by $t_0=a<t_1<\ldots <t_n=b$. Let $\xi_j\in I_{\Box j}=[t_{j-1},t_j]_{\Box}$ for
$j=1, 2, \ldots, n$. Then, $\sum_{j=1}^n f(\xi_j)\Delta g_j\leq
U_{\Box}(P,f,g)$ and $L_{\Box}(P,f,g)\leq \int_a^b f(t)\Box g(t)$. Hence,
\[ \left|\sum_{j=1}^nf(\xi_j)\Delta g_j-\int_a^b f(t)\Box g(t)\right|
\leq |U_{\Box}(P,f,g)-L_{\Box}(P,f,g)|<\eps\, . \]
\end{proof}
Theorem~\ref{integal_box} gives a comparison between the two kinds of Riemann--Stieljes integrals on time scales and
the Riemann--Stieljes integral on the real interval. In the particular case $g(t) = t$ one gets the result as stated in \cite{BPi}.
\begin{theorem}
\label{integal_box}
Let $a,b\in\T$ and $g:[a,b]\rightarrow\R$ be a strictly increasing function on $[a,b]$. Let $f:[a,b]\rightarrow\R$. Denote by $\left.f\right|_{\T}$ and $\left.g\right|_{\T}$ the restrictions of functions $f$ and $g$ to the time scale $\T$. Then,
\begin{enumerate}
  \item [i)] $\int_a^b \left.f\right|_{\T}(t)\Delta \left.g\right|_{\T}(t)\leq \int_a^b f(t)dg(t)\leq \int_a^b \left.f\right|_{\T}(t)\nabla \left.g\right|_{\T}(t)$
      if $f$ is strictly increasing on $[a,b]$;

  \item [ii)] $\int_a^b \left.f\right|_{\T}(t)\nabla \left.g\right|_{\T}(t)\leq \int_a^b f(t)dg(t)\leq
\int_a^b \left.f\right|_{\T}(t)\Delta \left.g\right|_{\T}(t)$
  if $f$ is strictly decreasing on $[a,b]$.
\end{enumerate}
\end{theorem}
\begin{proof}
We prove here the item i). The proof of ii) is similar. If $f$ is strictly increasing on $[a,b]$, then $\left.f\right|_{\T}$
is regulated on
$[a,b]_{\T}$ and all integrals exist for an increasing function $g$. Let $P\in \mathcal{P}\left([a,b]_{\T}\right)$ and $P=\{t_0,\ldots,
t_n\}$, where $a=t_0<t_1<\ldots<t_n=b$. Then, $\int_a^bf(t)dg(t)\leq
\sum_{j=1}^n f(t_{j})\Delta g_j= U_{\nabla}(P,f,g)$. Taking the infimum of the right-hand side over all partitions from
$\mathcal{P}\left([a,b]_{\T}\right)$  we get
$\int_a^bf(t)dg(t)\leq\overline{\int_a^b}\left.f\right|_{\T}(t)\nabla
\left.g\right|_{\T}(t)$.
Similarly, $\int_a^b f(t)dg(t)\geq \sum_{j=1}^n f(t_{j-1})\Delta g_j= L_{\Delta}(P,f,g)$. Taking now the supremum of the right-hand side we get that $\int_a^b f(t)dg(t)\geq\underline{\int_a^b}\left.f\right|_{\T}(t)\Delta
\left.g\right|_{\T}(t)$. For $\left.f\right|_{\T}\in
R_{\Box}(\left.g\right|_{\T},[a,b]_{\T})$, then also $\int_a^bf(t)dg(t)\leq
\int_a^b \left.f\right|_{\T}(t)\nabla \left.g\right|_{\T}(t)$ and $\int_a^b f(t)dg(t)\geq
\int_a^b \left.f\right|_{\T}(t)\Delta \left.g\right|_{\T}(t)$.
\end{proof}

\begin{corollary}
Let $a,b\in\T$, $I=[a,b]_{\T}$, $g: I \rightarrow\R$ be
a strictly increasing function, and $f: I \rightarrow\R$. Then,
\begin{enumerate}
  \item [i)] $\int_a^b f(t)\Delta g(t)\leq  \int_a^b f(t)\nabla g(t)$
if $f$ is strictly increasing on $I$;

  \item [ii)] $\int_a^b f(t)\nabla g(t)\leq  \int_a^b f(t)\Delta g(t)$ if  $f$
is strictly decreasing on $I$.
\end{enumerate}
\end{corollary}


\section{Algebraic properties of the Riemann-Stieltjes $\Box$--integral}

In this section we prove some algebraic properties of the
Riemann-Stieltjes integral on time scales. The properties are valid for an arbitrary time scale $\T$
with at least two points.
We define $\int_a^a f(t)\Box g(t)=0$ and $\int_a^b f(t)\Box g(t)=-\int_b^a f(t)\Box g(t)$ for $a>b$.
\begin{theorem}
Let $I=[a,b]_{\T}$, $a,b\in\T$. Every constant function
$f:\T\rightarrow\R,$ $f(t) \equiv c$,
is Stieltjes $\Box$--integrable with
respect to $g$ on $I$ and
 \[\int_a^bc\Box g(t)=c\left(g(b)-g(a)\right).\]
\end{theorem}
\begin{proof}
Let $P\in \mathcal{P}(I)$ and $P=\{t_0,\ldots,t_n\}$. Then,
\begin{equation*}
\begin{split}
 L_{\Box}(P,f,g) &=U_{\Box}P,f,g)=c\sum\limits_{j=1}^n\Delta g_j\\
 &= g(t_1)-g(a)+g(t_2)-g(t_1)+\cdots+g(b)-g(t_{n-1})\\
 &=g(b)-g(a).
\end{split}
\end{equation*}
Hence, $\underline{\int_a^b} f\Box g=\overline{\int_a^b} f\Box
g=c\left(g(b)-g(a)\right)$.
\end{proof}

\begin{theorem}
Let $t\in\T$ and  $f:\T\rightarrow\R.$ Then, $f$ is Riemann-Stieltjes
$\Delta$--integrable with respect to $g$ from $t$ to $\sigma(t)$ and
\begin{equation}
\label{eq:2}
 \int_t^{\sigma(t)}f(\tau) \Delta g(\tau)
 =f(t)\left(g^{\sigma}(t)-g(t)\right)\, ,
\end{equation}
where $g^{\sigma} = g \circ \sigma$.
Moreover, if $g$ is $\Delta$--differentiable at $t$, then
\begin{equation}
\int_t^{\sigma(t)}f(\tau) \Delta g(\tau)
=\mu(t)f(t)g^{\Delta}(t)\, . \label{eq:3}
\end{equation}
\end{theorem}
\begin{proof}
For $\sigma(t)=t$ both equations (\ref{eq:2}) and (\ref{eq:3})
hold. When $\sigma(t)>t$, only one partition for
$I=[t,\sigma(t)]_{\T}$ is possible. In that case
$P=\{t,\sigma(t)\}$ and we have one
set $I_{\Delta 1}=[t,\rho(\sigma(t))]_{\T}=\{t\}$. Hence, $m_{\Delta 1}=M_{\Delta 1}=f(t)$, and
\[\underline{\int_t^{\sigma(t)}} f\Delta g
=\overline{\int_t^{\sigma(t)}} f\Delta g
=U_{\Delta}(P,f,g)=L_{\Delta}(P,f,g)=f(t)(g(\sigma(t))-g(t)) \, .\]
As $\sigma(t)>t$, then
also $\mu(t)\neq 0$ and $g^{\sigma}(t)-g(t)=\mu(t)g^{\Delta}(t)$.
\end{proof}

\begin{theorem}
Let $t\in\T$ and  $f:\T\rightarrow\R.$ Then, $f$ is Riemann-Stieltjes
$\nabla$--integrable with respect to $g$ from $\rho(t)$ to $t$ and
\begin{equation}
\label{eq:25}
 \int_{\rho(t)}^{t}f(\tau) \nabla g(\tau)
 =f(t)\left(g(t)-g^{\rho}(t)\right)\, ,
\end{equation}
where $g^{\rho} = g \circ \rho$.
Moreover, if $g$ is $\nabla$--differentiable at $t$, then
\begin{equation}
\int_{\rho(t)}^{t}f(\tau) \nabla g(\tau)
=\nu(t)f(t)g^{\nabla}(t)\, . \label{eq:36}
\end{equation}
\end{theorem}
\begin{proof}
For $\rho(t)=t$ both equations (\ref{eq:25}) and (\ref{eq:36})
hold. When $t>\rho(t)$, only one partition for
$I=[\rho(t), t]_{\T}$ is possible. In that case
$P=\{\rho(t), t\}$ and we have one
set $I_{\nabla 1}=[\sigma(\rho(t)),t]_{\T}=\{t\}$. Hence, $m_{\nabla 1}=M_{\nabla 1}=f(t)$, and
\[\underline{\int_{\rho(t)}^{t}} f\nabla g
=\overline{\int_{\rho(t)}^{t}} f\nabla g
=U_{\nabla}(P,f,g)=L_{\nabla}(P,f,g)=f(t)(g(t)-g(\rho(t))) \, .\]
As $\rho(t)<t$, then
also $\nu(t)> 0$ and $g(t)-g^{\rho}(t)=\nu(t)g^{\nabla}(t)$.
\end{proof}

\begin{corollary}
Let $a,b\in\T$ and $a<b$.
\begin{enumerate}
\item [i)] If $\T=\R$, then a bounded function $f$ on $[a,b]_{\T}$ is
Riemann-Stieltjes $\Box$--integrable with respect to the
increasing function $g$ from $a$ to $b$ if and only if $f$ is
Riemann-Stieltjes integrable on $I$ in the classical sense.
Moreover, then  $\int_a^bf(t)\Box g(t)=\int_a^bf(t)dg(t)$,
where the integral on the right hand side is the classical
Riemann-Stieltjes integral (see, \textrm{e.g.}, \cite[Chapter~4]{B:garden}).

\item [ii)] If  $\T=\Z$, then each function $f:\Z\rightarrow\R$ is
Riemann-Stieltjes $\Box$--integrable from $a$ to $b$ with respect to an
arbitrarily increasing function $g:\Z\rightarrow\R$. Moreover,
$\int_a^bf(t)\Delta g(t)=\sum_{t=a}^{b-1}f(t)(g(t+1)-g(t))$,  $\int_a^bf(t)\nabla g(t)=\sum_{t=a+1}^{b}f(t)(g(t)-g(t-1))$.
\end{enumerate}
\end{corollary}
\begin{proof}
i) Notice that for $\T=\R$ Definition~\ref{RSdef} coincides
 with the classical definition of the Riemann-Stieltjes integral.
 Moreover, since for $\T=\R$ the $\Box$--differential coincides
 with the standard differential, it follows that
 \[\int_a^bf(t)\Box g(t)=\int_a^bf(t)g'(t)dt=\int_a^bf(t)dg(t) \, .\]

ii) Let $I=[a,b]_{\Z}$, $a<b$, $a$, $b\in\Z$. Consider the partition
  $P_{\delta} \in \mathcal{P}(I)$ given by $a=t_0<t_1<\ldots<t_n=b=a+n$,
  where $t_j=a+j$ for $j=0,\ldots, n$.
  Notice that $[t_{j-1},\rho(t_j)]_{\T}=\{t_{j-1}\}$ for each $j=1,\ldots, n$
  and $\Delta g(t_j)=g(t_j)-g(t_j-1)$. A direct calculation shows that
  $L_{\Delta}(P_{\delta}, f, g)=U_{\Delta}(P_{\delta},f,g)
  =\sum_{j=1}^nf(t_{j-1})\Delta g(t_j)
  =\sum_{t=a}^{b-1}f(t)(g(t+1)-g(t))$.
  Hence, $U_{\Delta}(P_{\delta},f,g)=L_{\Delta}(P_{\delta},f,g)$
  and Theorem~\ref{th.7} imply the desired formula for the $\Delta$--case. Similarly proof holds for the $\nabla$--case.
\end{proof}

\begin{theorem}(Linearity)
\label{lin}
Let functions $f_1$ and $f_2$ be Riemann-Stieltjes
$\Box$--integrable on the interval $[a,b]_{\T}$ with respect to  $g$, and $c$ be a constant. Then,
\begin{enumerate}
 \item [i)] $cf_1 \in \mathcal{R}_{\Box}(g,I)$
 and $\int_a^bcf_1 \Box g = c \int_a^b f_1\Box g$;

 \item [ii)] $f_1+f_2 \in \mathcal{R}_{\Box}(g,I)$
 and $\int_a^b\left(f_1+f_2\right)\Box g
 =\int_a^bf_1\Box g+\int_a^bf_2\Box g$.
\end{enumerate}
\end{theorem}
\begin{proof}
Let $I_j=(t_{j-1},t_j]_{\T}$, and
$p\in\mathcal{P}(I)$ be any partition of the interval
$I$. Denote $I_j=[t_{j-1},t_j]_{\Box}$ and
 \begin{eqnarray*}
  m_j^{f_1}=\inf\limits_{t\in I_{\Box  j}}f_1(t) \, , \quad
  \quad
  M_j^{f_1}=\sup\limits_{t\in I_{\Box  j}}f_1(t) \, ,\\
  m_j^{f_2}=\inf\limits_{t\in I_{\Box  j}}f_2(t)\, , \quad
  \quad
  M_j^{f_2}=\sup\limits_{t\in I_{\Box  j}}f_2(t) \, .
 \end{eqnarray*}
 \emph{i)} Let us notice that for $c=0$ function
 $cf_1 \in \mathcal{R}_{\Box}(g,I)$. Suppose that $c \neq 0$ and denote
 \[\mathfrak{m}_j=\inf\limits_{t\in I_{\Box  j}}(cf_1)(t)\quad
 \text{and}\quad
 \mathfrak{M}_j=\sup\limits_{t\in I_{\Box  j}}(cf_1)(t) \, .\]
 Then,
 \[\mathfrak{m}_j=cm_j^{f_1}\quad
 \text{and}\quad
 \mathfrak{M}_j=cM_j^{f_1}\quad
 \text{for}\;c>0;\]
 \[\mathfrak{m}_j=cM_j^{f_1}\quad
 \text{and}\quad
 \mathfrak{M}_j=cm_j^{f_1}\;\;\;\text{for}\;c<0\,.\]
 This implies that
 \[
 L_{\Box}(P,cf_1,g) =\left\{\begin{array}{ll}
					 cL_{\Box}(P,f_1,g), & \hbox{for\;$c>0$} \\
					 cU_{\Box}(P,f_1,g), & \hbox{for\;$c<0$}
\end{array}\right.\]
and
\[U_{\Box}(P,cf_1,g)=\left\{\begin{array}{ll}
  cU_{\Box}(P,f_1,g), & \hbox{for\;$c>0$} \\
  cL_{\Box}(P,f_1,g), & \hbox{for\;$c<0$.}
\end{array}\right.\]
Thus,
\begin{equation}
\label{war}
  U_{\Box}(P,cf_1,g)-L_{\Box}(P,cf_1,g)=|c|(U_{\Box}(P,f_1,g)-L_{\Box}(P,f_1,g))\,.
\end{equation}
Because $f_1\in\mathcal{R}_{\Box}(g,I)$, then for any $\varepsilon>0$ there
exists a partition $P$ such that
$U_{\Box}(P,cf_1,g)-L_{\Box}(P,cf_1,g)<\frac{\varepsilon}{|c|}$.
Together with (\ref{war}) this leads to
$U_{\Box}(P,cf_1,g)-L_{\Box}(P,cf_1,g)<\varepsilon$. Hence,
$cf_1 \in \mathcal{R}_{\Box}(g,I)$. Moreover,
 \[\underline{\int_a^b}cf_1 \Box g
 =\sup_{P\in \mathcal{P}(I)} L_{\Box}(P,cf_1,g) = c\sup _{P\in \mathcal{P}(I)}L_{\Box}(P,f_1,g)
 = c\underline{\int_a^b}f_1\Box g\]
 and
 \[\overline{\int_a^b}cf_1\Box g=\inf_{P\in \mathcal{P}(I)} U_{\Box}(P,cf_1,g) = c\inf_{P\in \mathcal{P}(I)} U_{\Box}(P,f_1,g)=
   c\overline{\int_a^b}f_1\Box g\,. \]
 Since $f_1 \in \mathcal{R}_{\Box}(g,I)$, then
 $\overline{\int_a^b}f_1\Box g=\underline{\int_a^b}f_1\Box g$,
 which proves the intended conclusion. \\

 \emph{ii)} Let
  \[\mathfrak{m}_j=\inf\limits_{t\in I_{\Box  j}}(f_1+f_2)(t)\quad
  \text{and}\quad
  \mathfrak{M}_j=\sup\limits_{t\in I_{\Box  j}} (f_1+f_2)(t) \, .\]
 Then, $\mathfrak{m}_j \geq m_j^{f_1}+m_j^{f_2}$
 and $\mathfrak{M}_j \leq M_j^{f_1}+M_j^{f_2}$.
 This implies
\begin{equation*}
\begin{gathered}
L_{\Box}(P,f_1,g)+L_{\Box}(P,f_2,g) \leq L_{\Box}(P,f_1+f_2,g)\\
U_{\Box}(P,f_1,g)-U_{\Box}(P,f_2,g) \leq U_{\Box}(P,f_1+f_2,g)\,.
\end{gathered}
\end{equation*}
 The integrability criterion implies the existence of partitions
 $P_1\in\mathcal{P}(I)$ and $P_2\in\mathcal{P}(I)$ such that the
 inequalities
 \begin{equation}
 \label{3.4}
  U_{\Box}(P_1,f,g)-L_{\Box}(P_1,f,g)<\frac{\varepsilon}{2}\quad
  \text{and}\quad
  U_{\Box}(P_2,f,g)-L_{\Box}(P_2,f,g)<\frac{\varepsilon}{2}
 \end{equation}
 also hold on their common refinement $P_{\varepsilon}=P_1\cup
 P_2$. From (\ref{3.4}) follows that
 $U_{\Box}(P_{\varepsilon},f_1+f_2,g)-L_{\Box}(P_{\varepsilon},f_1+f_2,g)<\varepsilon$.
 Hence, $f_1+f_2 \in \mathcal{R}_{\Box}(g,I)$. Additionally,
 \begin{multline*}
  \underline{\int_a^b}(f_1+f_2)\Box g=\sup_{P\in \mathcal{P}(I)} L_{\Box}(P,f_1+f_2,g) \\
   \leq \sup_{P\in \mathcal{P}(I)} L_{\Box}(P,f_1,g)+\sup_{P\in \mathcal{P}(I)} L_{\Box}(P,f_2,g)=
   \underline{\int_a^b}f_1\Box g+\underline{\int_a^b}f_2\Box g
 \end{multline*}
 and
 \begin{multline*}
  \overline{\int_a^b}(f_1+f_2)\Box g=\inf_{P\in \mathcal{P}(I)} U_{\Box}(P,f_1+f_2,g) \\
  \geq \inf_{P\in \mathcal{P}(I)} U_{\Box}(P,f_1,g)+\inf_{P\in \mathcal{P}(I)} U_{\Box}(P,f_2,g)=
   \overline{\int_a^b}f_1\Box g+\overline{\int_a^b}f_2\Box g \, ,
 \end{multline*}
 so that
 \[\int_a^b(f_1+f_2)\Box g \geq\overline{\int_a^b}f_1\Box g+\overline{\int_a^b}f_2\Box g \ \ \  \mbox{and} \ \ \
 \int_a^b(f_1+f_2)\Box g \geq
   \underline{\int_a^b}f_1\Box g+\underline{\int_a^b}f_2\Box g \, .\]
 Because $f_1 \in \mathcal{R}_{\Box}(g,I)$ and $f_2 \in \mathcal{R}_{\Box}(g,I)$,
 then $\overline{\int_a^b}f_1\Box g=\underline{\int_a^b}f_1\Box g$ and
 $\underline{\int_a^b}f_2\Box g=\overline{\int_a^b}f_2\Box g$.
\end{proof}

Similarly to the proof of Theorem~\ref{lin} item i),
one can show the following:

\begin{theorem}
Let $f\in \mathcal{R}_{\Box}(g_1,I)$ and $f\in \mathcal{R}_{\Box}(g_2, I)$,
where $g_1$ and $g_2$ are increasing functions on $[a,b]_{\T}$. Then,
$f \in \mathcal{R}_{\Box}(g_1+g_2,I)$ and $\int_a^bf\Box\left(g_1+ g_2\right)
= \int_a^bf\Box g_1+\int_a^bf\Box g_2$.
\end{theorem}

\begin{theorem}
Let $a$, $b$, $c\in\T$ with $a<b<c$. If $f$
is bounded on $[a,c]_{\T}$
and $g$ is monotonically increasing on $[a,c]_{\T}$, then
\begin{equation*}
\int_a^cf\Box g
= \int_a^bf\Box g+ \int_b^cf\Box g\, .
\end{equation*}
\end{theorem}
\begin{proof}
 There exist partitions $P_1 \in \mathcal{P}([a,b]_{\T})$ and $P_2 \in
 \mathcal{P}([b,c]_{\T})$ such that for $\varepsilon>0$
 \[U_{\Box}(P_1,f,g) - L_{\Box}(P_1,f,g)<\frac{\varepsilon}{2}\quad
 \text{and}\quad
 U_{\Box}(P_2,f,g) -
 L_{\Box}(P_2,f,g)<\frac{\varepsilon}{2} \, .\]
Then, there is a partition $P \in \mathcal{P}([a,c]_{\T})$ such that
 $U_{\Box}(P,f,g)=U_{\Box}(P_1,f,g)+U_{\Box}(P_2,f,g)$ and
 $L_{\Box}(P,f,g)=L_{\Box}(P_1,f,g)+L_{\Box}(P_2,f,g)$ so that
 $U_{\Box}(P,f,g)-L_{\Box}(P,f,g)<\varepsilon$. Hence, $f$ is $\Box$--integrable
 with respect to $g$ and
 \begin{multline*}
 \int_a^c f \Box g \leq U_{\Box}(P_1,f,g)+U_{\Box}(P_2,f,g) \\
 \leq  L_{\Box}(P_1,f,g)+L_{\Box}(P_2,f,g)+\varepsilon
 \leq \int_a^b f \Box g +\int_b^c f \Box g +\varepsilon \, .
 \end{multline*}
Similarly,
 \[\int_a^c f \Box g \geq \int_a^b f \Box g +\int_b^c f \Box g -\varepsilon \, .\]
\end{proof}

\begin{theorem}(Integration by substitution)
Let $\tilde{\T}$ be a time scale.
Assume that $\varphi:\ti{\T} \rar \R$ is a strictly increasing
continuous function and $\T=\varphi(\ti{\T})$ is a time scale.
Moreover, $\varphi$ maps an interval $[A, B]_{\tilde{\T}}$ onto $[a,b]_{\T}$.
Let $g$ be monotonically increasing on $[a, b]_{\T}$ and
$f\in \mathcal{R}_{\Box}(g,[a,b]_{\T})$. Then,
$f\circ\varphi\in \mathcal{R}_{\Box}(g\circ\varphi, [A,B]_{\tilde{\T}})$ and
\begin{equation*}
\int_a^b f(t) \Box g(t)=\int_A^B f(\varphi(s))\Box g(\varphi(s)).
\end{equation*}
\end{theorem}
\begin{proof}
Since $\varphi$ is a strictly increasing continuous function,
there is a one-to-one correspondence between the partition
 $Q=\{s_0, s_1,\ldots, s_n\}\subset [A, B]_{\tilde{\T}}$,
 where $A=s_0<s_1< \ldots <s_n=B$,
 and the partition
 $P=\{t_0, t_1, \ldots, t_n\}\subset [a,b]_{\T}$,
 where $a=t_0=\varphi(s_0)<t_1=\varphi(s_1)< \ldots <t_n=\varphi(s_n)=b$.
 Since $f([t_{j-1},t_j]_{\T})=f([\varphi(s_{j-1},\varphi(t_j)]_{\tilde{\T}})$
 for each $j$, then $L_{\Box}(P,f,g)=L_{\Box}(Q,f\circ\varphi, g\circ\varphi)$
 and $U_{\Box}(P,f,g)=U_{\Box}(Q,f\circ\varphi, g\circ \varphi)$.
 The result follows from Theorem~\ref{dor:criterion}.
\end{proof}


\section{From the Riemann-Stieltjes to Riemann $\Box$--integral}

Theorem~\ref{dor:transition} below establishes a relation between
the Riemann-Stieltjes $\Box$--integral and the Riemann $\Box$--integral for a function $g$ being $\Box$--differentiable on the interval $I$ of integration.
We begin by noting that one can easily reformulate to nabla version, with the interval opened from the left side,
the delta mean value theorem \cite[Theorem~1.14]{Bh1}:

\begin{theorem}
\label{thm:Th1.14:nabla}
Let $f$ be a continuous function on $[a, b]_{\T}$
that is $\nabla$ -- differentiable on $(a, b]$. Then there exist $\xi$, $\tau\in (a,b]$ such that
\[f^{\nabla}(\xi)\leq \frac{f(b)-f(a)}{b-a}\leq f^{\nabla}(\tau)\, .\]
\end{theorem}

In Corollary~\ref{box_mean} we write together
\cite[Theorem~1.14]{Bh1} and Theorem~\ref{thm:Th1.14:nabla}
with our "$\Box$"--notation:

\begin{corollary}
\label{box_mean}
Let $f$ be a continuous function on $[a, b]_{\T}$
that is $\Box$ -- differentiable on $[a, b]_{\Box}$. Then there exist $\xi$, $\tau\in[a, b]_{\Box}$ such that
\[f^{\Box}(\xi)\leq \frac{f(b)-f(a)}{b-a}\leq f^{\Box}(\tau)\, .\]
\end{corollary}

\begin{theorem}
\label{dor:transition}
Let $I=[a,b]_{\T}, a,b\in\T$. Suppose that $g$ is an increasing
function such that $g^{\Box}$ is continuous on $(a,b)_{\T}$
and $f$ is a real bounded function on $I$.
Then, $f\in\mathcal{R}_{\Box}(g, I)$ if and only if
$fg^{\Box}\in \mathcal{R}_{\Box}(g,I)$. Moreover,
\begin{equation*}
 \int_a^bf(t)\Box g(t)=\int_a^bf(t)g^{\Box}(t)\Box t \, .
\end{equation*}
\end{theorem}
\begin{proof}
 Let $\eps>0$. Since $g^{\Box}$ is continuous, then $g^{\Box}\in \mathcal{R}_{\Box}(t,I)$,
\textrm{i.e.}, $g^{\Box}$ is Riemann--Stieltjes $\Box$--integrable on $I$. Hence, there
 exists a partition $P=\{a=t_0, t_1,\ldots, t_n=b\}\in
 \mathcal{P}(I)$ such that $U_{\Box}(P,g^{\Delta},t)- L_{\Box}(P,g^{\Delta},t)<\frac{\eps}{M}$,
 where $M=\sup_{t\in I}|f(t)|$.
From the $\Box$--version of the mean value theorem on time scales (see Corollary~\ref{box_mean}),
 for each $j=1,\ldots,n$ there are
 $\xi_j, \tau_j\in I_j=[t_{j-1},t_j]_{\Box}$ such that
 \[g^{\Box}(\tau_j)\Delta t_j\leq \Delta g_j\leq  g^{\Box}(\xi_j)\Delta t_j \ .\]
 From Theorem~\ref{dor:3} part ii), for any $s_j\in I_j$, $j=1, \ldots, n$,
 \[\sum_{j=1}^n \left|g^{\Box}(\xi_j)-g^{\Box}(s_j)\right|\Delta t_j<\eps\, .\]
 Then, $\sum_{j=1}^n f(s_j)\Delta g_j\leq \sum_{j=1}^n f(s_j)g^{\Box}(\xi_j)\Delta t_j$ and
 \[\left|\sum_{j=1}^n f(s_j)\Delta g_j -\sum_{j=1}^n f(s_j)g^{\Box}(s_j)\Delta t_j\right|\leq
 \left|\sum_{j=1}^n f(s_j)g^{\Box}(\xi_j)\Delta t_j-\sum_{j=1}^n f(s_j)g^{\Box}(s_j)\Delta t_j\right| \, .\]
 Hence,
 \begin{equation}
 \label{dor1}
  \left|\sum_{j=1}^n f(s_j)\Delta g_j -\sum_{j=1}^n f(s_j)g^{\Box}(s_j)\Delta t_j\right|
  \leq M\left|\sum_{j=1}^n \left(g^{\Box}(\xi_j)-g^{\Box}(s_j)\right) \Delta t_j\right|<\eps
 \end{equation}
 for any $s_j\in I_j$. Thus, $\sum_{j=1}^n f(s_j)\Delta g_j\leq U_{\Box}(P, fg^{\Box},t)+\eps$
 and $U_{\Box}(P,f,g)\leq U_{\Box}(P,fg^{\Box},t)+\eps$. Inequality (\ref{dor1}) implies that
 $U_{\Box}(P,fg^{\Box},t)\leq U_{\Box}(P,f,g)+\eps$. Thus,
 \[\left|U_{\Box}(P,f,g)-U_{\Box}(P,fg^{\Box},t)\right|\leq \eps \,.\]  Moreover,
 \[\left|\overline{\int_a^b}f(t)\Box g(t)-\overline{\int_a^b}f(t)
 g^{\Box}(t)\Box t \right|\leq \eps\, .\]
 Since $\eps$ is arbitrary, we conclude that
 \[\overline{\int_a^b}f(t)\Box g(t)=\overline{\int_a^b}f(t) g^{\Box}(t)\Box t.\]
 In a similar way one prove that
 $\underline{\int_a^b}f(t)\Box g(t)
 =\underline{\int_a^b}f(t) g^{\Box}(t)\Box t$.
 \end{proof}

\begin{theorem}(Delta integration by parts)
Let $I=[a,b]_{\T}, a,b\in\T$. Suppose that $g$ is an increasing
function such that $g^{\Delta}$ is continuous on $(a,b)_{\T}$
and $f$ is a real bounded function on $I$. Then,
\begin{equation*}
\int_a^bf\Delta g=\left[fg\right]_a^b-\int_a^b g^{\sigma}\Delta f \, ,
\end{equation*}
 where, as usual, $g^{\sigma}$ means $g \circ \sigma$.
\end{theorem}
\begin{proof}
 Theorem~\ref{dor:transition} and integration by parts for
 the Riemann $\Delta$--integral on time scales (see \cite{Bh1}) imply
 that $\int_a^bf(t)\Delta g(t)=\int_a^bf(t)g^{\Delta}(t)\Delta t$
 and
$\int_a^bg^{\sigma}f^{\Delta}\Delta t +\int_a^bf\Delta t=[fg]_a^b$.
 Hence, $\int_a^bf\Delta g=\left[fg\right]_a^b-\int_a^b g^{\sigma}\Delta f$.
\end{proof}

\begin{theorem}(Nabla integration by parts)
Let $I=[a,b]_{\T}, a,b\in\T$. Suppose that $g$ is an increasing
function such that $g^{\nabla}$ is continuous on $(a,b)_{\T}$
and $f$ is a real bounded function on $I$. Then,
\begin{equation*}
\int_a^bf\nabla g=\left[fg\right]_a^b-\int_a^b g^{\rho}\nabla f \, ,
\end{equation*}
 where  $g^{\rho}=g \circ \rho$.
\end{theorem}
\begin{proof}
 Theorem~\ref{dor:transition} and integration by parts for
 the Riemann $\nabla$--integral on time scales  implies
 that $\int_a^bf(t)\nabla g(t)=\int_a^bf(t)g^{\nabla}(t)\nabla t$
 and
$\int_a^bg^{\rho}f^{\nabla}\nabla t +\int_a^bf\nabla t=[fg]_a^b$.
 Hence, $\int_a^bf\nabla g=\left[fg\right]_a^b-\int_a^b g^{\rho}\nabla f$.
\end{proof}


\section{Conclusion and Future Perspectives}
\label{sec:ConcFuture}

This article is about the concept of
Riemann-Stieltjes delta integration on time scales.
The results of the paper may be used, \textrm{e.g.},
to generalize the $\T = \R$ inequalities
proved in \cite{D:2004,D:2002}
to a general time scale $\T$. Then,
as the particular case $g(t)=t$, one would
obtain the previous inequalities
proved for the Riemann integral
on time scales \cite{R:R:D:JIA,R:D:JMS:2009}.

Another interesting line of research is to
investigate the possibility of extending all
previous notions of integration on time scales
by putting together our present results with
the Henstock-Kurzweil integrals introduced
by Peterson and Thompson in \cite{P:2006,Thomson}.
Such Henstock-Kurzweil-Stieltjes integrals on time scales
are under study and will be addressed elsewhere.


\begin{acknowledgement}
The first author was supported by Bia\l ystok Technical University through W/WI/7/07 grant; the second and third authors by the Portuguese Foundation for Science and Technology
(Paw\l uszewicz through the program ``1 000 Doutorados
para Institui\c{c}\~{o}es Cient\'{\i}ficas Portuguesas'';
Torres through the R\&D unit CEOC of the University
of Aveiro, cofinanced by the European Community fund
FEDER/POCI 2010).
\end{acknowledgement}



\end{document}